 \documentclass[final,1p,times]{elsarticle}

\usepackage{amssymb}
\usepackage{amsmath}

\newtheorem{definition}{Definition}[section]
\newtheorem{theorem}{Theorem}[section]
\newtheorem{corollary}[theorem]{Corollary} 

\newtheorem{lemma}{Lemma}[section]

\numberwithin{equation}{section}
\journal{   }
\begin{document}
	
	\begin{frontmatter}
		
		
		
		\title{Positive definiteness of fourth order three dimensional symmetric tensors} 
		
		
		\author{Yishenhg Song}
		
		\address{School of Mathematical Sciences, Chongqing Normal University, Chongqing 401331 P.R. China.\\
			Email: yisheng.song@cqnu.edu.cn}
		
		
		\begin{abstract}
			For a 4th order 3-dimensional symmetric tensor with its entries  $1$ or $-1$, we show the analytic  sufficient and necessary  conditions  of  its positive definiteness. By applying these conclusions, several  strict inequalities is bulit for ternary quartic homogeneous polynomials.
		\end{abstract}
		
		
		
		\begin{keyword}
			Positive definiteness, Fourth order tensors, Homogeneous polynomial.
			
			
		\end{keyword}
		
	\end{frontmatter}
\footnotetext{The author’s work was supported by the National Natural Science Foundation of P.R. China (Grant No.12171064), by The team project of innovation leading talent in chongqing (No.CQYC20210309536) and by the Foundation of Chongqing Normal univer- sity (20XLB009).}

\section{Introduction}

One of the most direct applications of positive definite tensors is to verify   the  vacuum stability of the Higgs scalar potential  model  \cite{K2016,SQ2024}.  Qi \cite{Q2005} first used  the consept of  positive definiteness for a symmetric tensor when the order is  even integer.  

\begin{definition} Let $\mathcal{T}=(t_{i_1i_2\cdots i_m})$  be an
	$m$th order $n$ dimensional symmetric tensor.  $\mathcal{T}$ is called 
	\begin{itemize}
		\item[(i)] {\bf positive semi-definite} (\cite{Q2005}) if $m$ is  an even number and  in the Euclidean space $ \mathbb{R}^n$, its associated Homogeneous polynomial $$\mathcal{T}x^m=\sum\limits_{i_1,i_2,\cdots,i_m=1}^nt_{i_1i_2\cdots
			i_m}x_{i_1}x_{i_2}\cdots
		x_{i_m}\geq0;$$ 
		\item[(ii)] {\bf positive definite} (\cite{Q2005}) if $m$ is  an even number and  $\mathcal{T}x^m>0$ for all $x\in \mathbb{R}^n\setminus\{0\}$.
\end{itemize} 
\end{definition}

Clearly, a positive semi-definite tensor coincides with a positive semi-definite matrix if $m=2$. It is well-known that  Sylvester's Criterion can efficiently identify the positive (semi-)definiteness of a matrix.  The positive definiteness of a 4th order 2 dimensional symmetric tensor, (or positivity condition of a quartic univariate polynomial) may trace  back to ones of  Refs. Rees \cite{R1922} ,   Lazard \cite{L1988} Gadem-Li \cite{GL1964},  Ku \cite{K1965} and Jury-Mansour \cite{JM1981}.  Untill to 2005, Wang-Qi \cite{WQ2005} improved their proof, and perfectly gave  analytic necessary and sufficient conditions. However, the above result depends on the discriminant of such a quartic polynomial. Hasan-Hasa \cite{HH1996} claimed that a necessary and sufficient condition of positive definiteness  was proved without the discriminant. However, there is a problem in their argumentations. In 1998, Fu \cite{F1998} pointed out that Hasan-Hasan's results are sufficient only.  Recently, Guo\cite{G2021} showed a new necessary and sufficient condition without the discriminant. Very recently, Qi-Song-Zhang\cite{QSZ2022} gave a new necessary and sufficient condition other than the above results.  For more  detail about  applications of these results,  see Song-Qi \cite{SQ2021} also.  \\

In 2005, Qi \cite{Q2005} gave that the sign of all H-(Z-)eigenvalue of  a even order symmetric tensor can verify the positive definitiveness of such a higher order tensor.  Subsequently, Ni-Qi-Wang \cite{NQW2008}  provided  a method of computing the smallest eigenvalue for checking positive definiteness of a 4th order 3 dimensional tensor. Ng-Qi-Zhou \cite{NQZ2009} presented  an algorithm of the largest eigenvalue of a nonnegative tensor. For a 4th order 3 dimensional symmetric tensor, Song \cite{S2021} proved several  sufficient conditions of its positive definiteness. Until now, an analytic necessary and sufficient condition   has not  been  found  for positive (semi-)definiteness for a 4th order 3 dimensional symmetric tensor.

In this paper,  we mainly dicuss analytic  necessary and sufficient conditions of positive definiteness of a class of 4th order 3-dimensional symmetric tensors (Theorem \ref{thm:31}).  Furthermore,   several  strict inequalities of ternary quartic homogeneous polynomial  (Corollary \ref{cor:32}) are built.

\section{Copositivity of 4th order 2-dimensional symmetric tensors}

Let $\mathcal{T}=(t_{ijkl})$ be a 4th-order 2-dimensional symmetric tensor. Then for $x=(x_1,x_2)^\top,$
\begin{equation}\label{eq:f}
	Tx^4=t_{1111} x_1^4+4t_{1112} x_1^3x_2+6t_{1122} x_1^2x_2^2+4t_{1222}x_1x_2^3+t_{2222}x_2^4.
\end{equation}
Let $$\begin{aligned}
	\Delta	=&4\times 12^3(t_{1111}t_{2222}-4t_{1112}t_{1222}+3t_{1122}^2)^3\\
	&-72^2\times 6^2(t_{1111}t_{1122}t_{2222}+ 2t_{1112}t_{1122}t_{1222}- t_{1122}^3-t_{1111}t_{1222}^2- t_{1112}^2t_{2222})^2\\
	=& 4\times 12^3(I^3-27J^2),
\end{aligned}$$
where $$\begin{aligned}I=&t_{1111}t_{2222}-4t_{1112}t_{1222}+3t_{1122}^2,\\
	J=&t_{1111}t_{1122}t_{2222}+2t_{1112}t_{1122}t_{1222}-t_{1122}^3-t_{1111}t_{1222}^2-t_{1112}^2t_{2222}.
\end{aligned}$$
and hence, the sign of $\Delta$ is the same as one of  $(I^3-27J^2)$. Ulrich-Watson \cite{UW1994}  presented the analytic  conditions of the nonnegativity of a quartic and univariate polynomial   in $\mathbb{R}_+$.   Qi-Song-Zhang \cite{QSZ2022} also gave  the nonnegativity and positivity of a quartic and univariate polynomial   in $\mathbb{R}$, which means  the positive (semi-)definitiveness of 4th order 2-dimensional  tensor \cite{SQ2024}. 
\begin{lemma}[\cite{SQ2024,QSZ2022}]\label{lem:21} A 4th-order 2-dimensional symmetric tensor $\mathcal{T}=(t_{ijkl})$ is positive definite if and only if
	$$\begin{cases}
		I^3-27J^2=0,\ \ t_{1112}\sqrt{t_{2222}}=t_{1222}\sqrt{t_{1111}},\\
		2t_{1112}^2+t_{1111}\sqrt{t_{1111}t_{2222}}=3t_{1111}t_{1122}<3t_{1111}\sqrt{t_{1111}t_{2222}};\\
		I^3-27J^2>0,\ \ 
		|t_{1112}\sqrt{t_{2222}}-t_{1222}\sqrt{t_{1111}}|\leq \sqrt{6t_{1111}t_{1122}t_{2222}+2\sqrt{(t_{1111}t_{2222})^3}},\\
		(i) \ -\sqrt{t_{1111}t_{2222}}< 3t_{1221}\leq 3\sqrt{t_{1111}v_{2222}};\\
		(ii)\  t_{1221} >\sqrt{t_{1111}t_{2222}}\ \mbox{ and } \\ |t_{1112}\sqrt{t_{2222}}+t_{1222}\sqrt{t_{1111}}|\leq \sqrt{6t_{1111}t_{1122}t_{2222}-2\sqrt{(t_{1111}t_{2222})^3}}.
	\end{cases} \leqno{(\textbf{I})}$$
	A 4th-order 2-dimensional symmetric tensor $\mathcal{T}=(t_{ijkl})$ is positive semidefinite if and only if
	$$\begin{cases}
		I^3-27J^2\ge0,\ \ 
		|t_{1112}\sqrt{t_{2222}}-t_{1222}\sqrt{t_{1111}}|\leq \sqrt{6t_{1111}t_{1122}t_{2222}+2\sqrt{(t_{1111}t_{2222})^3}},\\
		(i) \ -\sqrt{t_{1111}t_{2222}}\leq 3t_{1122}\leq 3\sqrt{t_{1111}t_{2222}};\\
		(ii) 	\ t_{1122} >\sqrt{t_{1111}t_{2222}}\mbox{ and }\\
		|t_{1112}\sqrt{t_{2222}}+t_{1222}\sqrt{t_{1111}}|\leq \sqrt{6t_{1111}t_{1122}t_{2222}-2\sqrt{(t_{1111}t_{2222})^3}}.
	\end{cases}\leqno{(\textbf{II})}$$
\end{lemma}

\begin{lemma}\label{lem:22} Let $\mathcal{T}=(t_{ijkl})$  be a 4th-order 2-dimensional symmetric tensor with its entires $|t_{ijkl}|=1$ and $t_{1111}=t_{2222}=1.$  Then \begin{itemize}
		\item[(i)]  $\mathcal{T}$ is positive semidefinite if and only if
		$t_{1122}=1;$
		\item[(ii)]  $\mathcal{T}$ is positive definite if and only if
		$t_{1122}=1\mbox{ and } t_{1112}t_{1222}=-1.$
	\end{itemize}
\end{lemma}
{\bf Proof.} (i) It follows from Lemma \ref{lem:21} (II) that $\mathcal{T}$ is positive semidefinite if and only if
	$$I^3-27J^2\ge0,\ \  |t_{1112}-t_{1222}|\leq\sqrt{6t_{1122}+2}\mbox{ and }-1\leq 3t_{1122} \leq3. $$
	Since  $|t_{ijkl}|=1$, then which means 
	$t_{1122}=1$ and either $t_{1112}t_{1222}=1$,  $$ I^3-27J^2=(1-4+3)^{3}-27(1+2-1-1-1)^{2}=0,$$$$ |t_{1112}-t_{1222}|=0<\sqrt{6t_{1122}+2}=\sqrt8;$$ or  $t_{1112}t_{1222}=-1,$
	$$ I^3-27J^2=(1+4+3)^{3}-27(1-2-1-1-1)^{2}>0,$$$$ |t_{1112}-t_{1222}|=2<\sqrt{6t_{1122}+2}=\sqrt8.$$
	So $\mathcal{T}$ is positive semidefinite if and only if
	$t_{1122}=1.$
	
	(ii) It follows from Lemma \ref{lem:21} (I) that $\mathcal{T}$ is positive definite if and only if
	$$\aligned
	I^3-27J^2=&0,\ \  t_{1112}=t_{1222}, \ \ 2t_{1112}^2+1=3t_{1122}<3;\\
	I^3-27J^2>&0,\ \  |t_{1112}-t_{1222}|\leq\sqrt{6t_{1122}+2}\mbox{ and }-1< 3t_{1122} \leq3.
	\endaligned$$
	Since $3=2t_{1112}^2+1=3t_{1122}<3$ can't hold, then $\mathcal{T}$ is positive definite if and only if
	$t_{1122}=1\mbox{ and } t_{1112}t_{1222}=-1.$
	This  completes the proof.

\section{Positive definiteness of 4th order 3-dimensional symmetric tensors}
\begin{theorem}\label{thm:31} \em Let $\mathcal{T}=(t_{ijkl})$ be a 4th-order $3$-dimensional symmetric tensor with its entries, \begin{center}
		$|t_{iiij}|=|t_{iijk}|=t_{iiii}=1$ and $t_{ijjj}t_{iiij}=-1$ for all $i,j,k\in\{1,2,3\}$, $i\ne j$, $i\ne k$, $j\ne k.$
	\end{center}
	\begin{itemize}
		\item[(i)]  If $t_{iijj} =\dfrac{11}6$ for all $i,j\in\{1,2,3\}$ and $i\ne j$, then $\mathcal{T}$ is  positive semidefinite if and only if
		$$
		t_{1222}=t_{2333}=t_{1113}, t_{1112}=t_{1333}=t_{2223}, t_{1123}=t_{1223}=t_{1233}=-1.
		\leqno{(\textbf{III})}$$
		\item[(ii)]    If $t_{iijj} =2$ for all $i,j\in\{1,2,3\}$ and $i\ne j$, then $\mathcal{T}$ is  positive definite if and only if the above condition (\textbf{III}) holds. 
		\item[(iii)]  If $t_{iijj} =2.5$ for all $i,j\in\{1,2,3\}$ and $i\ne j$, then $\mathcal{T}$ is  positive definite if and only if
		$$\begin{cases}
			t_{1123}=t_{1223}=t_{1233}=1; \mbox{ or}\\
			t_{1123}=t_{1223}=t_{1233}=-1, t_{1222}=t_{2333}=t_{1113}, t_{1112}=t_{1333}=t_{2223}; \mbox{ or}\\
			\mbox{two of } \{t_{1123}, t_{1223}, t_{1233}\} \mbox{ are  }-1.
		\end{cases} \leqno{(\textbf{IV})}$$
		\item[(iv)]  If $t_{iijj} \ge\dfrac83$ for all $i,j\in\{1,2,3\}$ and $i\ne j$, then $\mathcal{T}$ is  positive definite. 
	\end{itemize} 
\end{theorem}
{\bf Proof.}  
$\mathcal{T}x^4 $ may be rewritten as follows, $$\aligned 
\mathcal{T}x^4=&(x_1+x_2+x_3)^4-8(x_1^3x_2+x_1x_3^3+x_2^3x_3)\\&+6(t_{1122}-1)x_1^2x_2^2+6(t_{1133}-1)x_1^2x_3^2+6(t_{2233}-1)x_2^2x_3^2\\
&+12(t_{1123}-1)x_1^2x_2x_3+12(t_{1223}-1)x_1x_2^2x_3+12(t_{1233}-1)x_1x_2x_3^2\\
=&(x_1+x_2-x_3)^4+8(x_1^3x_3+x_2x_3^3-x_1^3x_2)\\&+6(t_{1122}-1)x_1^2x_2^2+6(t_{1133}-1)x_1^2x_3^2+6(t_{2233}-1)x_2^2x_3^2\\
&+12(t_{1123}+1)x_1^2x_2x_3+12(t_{1223}+1)x_1x_2^2x_3+12(t_{1233}-1)x_1x_2x_3^2\\
=&(x_1-x_2+x_3)^4+8(x_1x_2^3+x_2x_3^3-x_1^3x_3)\\&+6(t_{1122}-1)x_1^2x_2^2+6(t_{1133}-1)x_1^2x_3^2+6(t_{2233}-1)x_2^2x_3^2\\
&+12(t_{1123}+1)x_1^2x_2x_3+12(t_{1223}-1)x_1x_2^2x_3+12(t_{1233}+1)x_1x_2x_3^2\\
=&(x_2+x_3-x_1)^4+8(x_1^3x_3+x_1x_2^3-x_2^3x_3)\\&+6(t_{1122}-1)x_1^2x_2^2+6(t_{1133}-1)x_1^2x_3^2+6(t_{2233}-1)x_2^2x_3^2\\
&+12(t_{1123}-1)x_1^2x_2x_3+12(t_{1223}+1)x_1x_2^2x_3+12(t_{1233}+1)x_1x_2x_3^2.
\endaligned$$
(i)   Necessity. 	Suppose the  conditions (\textbf{III}) can't hold when $\mathcal{T}$ is positive semidefinite,  then there may be  four cases. 

Case 1.  There is two  $-1$  in $\{t_{1123}, t_{1223}, t_{1233}\}$.  We might as well take $t_{1123}=t_{1223}=-1$ and $t_{1233}=1$. Without loss the generality, let $t_{1112}=t_{2333}=t_{1113}=1$ and $ t_{1222}=t_{1333}=t_{2223}=-1$.  Take $x=(\frac15,-\frac15,1)^\top$. Then we have
$$\aligned 
\mathcal{T}x^4=&(x_1+x_2-x_3)^4+8(x_1^3x_3+x_2x_3^3-x_1x_2^3)+5(x_1^2x_2^2+x_1^2x_3^2+x_2^2x_3^2)\\
=&1+8\left(\dfrac1{5^3}-\dfrac15+\dfrac1{5^4}\right)+5\left(\dfrac1{5^4}+\dfrac1{5^2}+\dfrac1{5^2}\right)=-\dfrac{72}{625}<0;
\endaligned$$

Case 2.	 There is only  one  $-1$  in $\{t_{1123}, t_{1223}, t_{1233}\}$. We might as well take $t_{1123}=t_{1233}=1$ and $t_{1223}=-1 $ and  $t_{1112}=t_{2333}=t_{1113}=1$ and $ t_{1222}=t_{1333}=t_{2223}=-1$.   For $x=(\frac12,-\frac12,1)^\top$, we have
$$\aligned 
\mathcal{T}x^4=&(x_1+x_2+x_3)^4-8(x_1x_2^3+x_1x_3^3+x_2^3x_3)+5(x_1^2x_2^2+x_1^2x_3^2+x_2^2x_3^2)-24x_1x_2^2x_3\\\
=&1-8\left(-\dfrac1{2^4}+\dfrac12-\dfrac1{2^3}\right)+5\left(\dfrac1{2^4}+\dfrac1{2^2}+\dfrac1{2^2}\right)-3=-\dfrac{27}{16}<0.
\endaligned$$

Case 3. $t_{1123}=t_{1223}=t_{1233}=1$ and  $t_{1112}=t_{2333}=t_{1113}=1$ and $ t_{1222}=t_{1333}=t_{2223}=-1$. Take $x=(\frac15,-\frac15,1)^\top$. Then we have $$ \aligned 
\mathcal{T}x^4=&(x_1+x_2+x_3)^4-8(x_1x_2^3+x_1x_3^3+x_2^3x_3)+5(x_1^2x_2^2+x_1^2x_3^2+x_2^2x_3^2)\\
=&1-8\left(-\dfrac1{5^4}+\dfrac15-\dfrac1{5^3}\right)+5\left(\dfrac1{5^4}+\dfrac1{5^2}+\dfrac1{5^2}\right)=-\dfrac{72}{625}<0.
\endaligned$$
The above three cases imply that the equality, $t_{1123}=t_{1223}=t_{1233}=-1,$ is necessary.

Case 4. $t_{1123}=t_{1223}=t_{1233}=-1,$ but $ t_{1222}=t_{2333}=t_{1113}$ and $ t_{1112}=t_{1333}=t_{2223}$ can't hold. Without loss the generality, let $t_{1112}=t_{2333}=t_{1113}=1$ and $ t_{1222}=t_{1333}=t_{2223}=-1$. Take $x=(-1,-3,-1)^\top$. Then we have $$ \aligned 
\mathcal{T}x^4=&(x_1+x_2-x_3)^4+8(x_1^3x_3+x_2x_3^3-x_1x_2^3)+5(x_1^2x_2^2+x_1^2x_3^2+x_2^2x_3^2)-24x_1x_2x_3^2\\
=&(-1-3+1)^4-8(1+3-3^3)+5(9+1+9)-24\times3=-80<0.
\endaligned$$
This is a contradiction to  the positive semidefiniteness of $\mathcal{T}$, and hence, the  conditions (\textbf{III})  are  necessary.

Sufficiency.  $t_{1123}=t_{1223}=t_{1233}= -1$ and $ t_{1222}=t_{2333}=t_{1113}$ and $ t_{1112}=t_{1333}=t_{2223}$.  Without loss the generality, let $t_{1222}=t_{2333}=t_{1113}=1$ and $ t_{1112}=t_{1333}=t_{2223}=-1$. Rewriting  $\mathcal{T}x^4 $ as follow, 
$$ 
\mathcal{T}x^4=(x_1+x_2-x_3)^4+8(x_1^3x_3+x_2x_3^3-x_1^3x_2)+5(x_1^2x_2^2+x_1^2x_3^2+x_2^2x_3^2)-24x_1x_2x_3^2.
$$
Solve the constrained optimization problem:
$$\aligned
\min \ &\ \mathcal{T}x^4 \\
\mbox{ s. t. } & x_1^2+x_2^2+x_3^2=1.
\endaligned$$
Then the minimum value   is $0$ at a point $\left(\frac1{\sqrt{3}},\frac1{\sqrt{3}},\frac1{\sqrt{3}}\right)$ or $\left(-\frac1{\sqrt{3}},-\frac1{\sqrt{3}},-\frac1{\sqrt{3}}\right)$, and  hence, $\mathcal{T}x^4\ge0$. That is, $\mathcal{T}$ is positive semidefinite.

(ii)   Necessity. 	Suppose the  conditions (\textbf{III}) can't hold when $\mathcal{T}$ is positive definite,  then there may be  four cases. 

Case 1.  There is two  $-1$  in $\{t_{1123}, t_{1223}, t_{1233}\}$.  We might as well take $t_{1123}=t_{1223}=-1$ and $t_{1233}=1 $ and  $t_{1112}=t_{2333}=t_{1113}=1$ and $ t_{1222}=t_{1333}=t_{2223}=-1$.  Take $x=(\frac15,-\frac15,1)^\top$. Then we have
$$\aligned 
\mathcal{T}x^4=&(x_1+x_2-x_3)^4+8(x_1^3x_3+x_2x_3^3-x_1x_2^3)+6(x_1^2x_2^2+x_1^2x_3^2+x_2^2x_3^2)\\
=&1+8\left(\dfrac1{5^3}-\dfrac15+\dfrac1{5^4}\right)+6\left(\dfrac1{5^4}+\dfrac1{5^2}+\dfrac1{5^2}\right)=-\dfrac{21}{625}<0;
\endaligned$$

Case 2.	 There is only  one  $-1$  in $\{t_{1123}, t_{1223}, t_{1233}\}$. We might as well take $t_{1123}=t_{1233}=1$ and $t_{1223}=-1 $ and  $t_{1112}=t_{2333}=t_{1113}=1$ and $ t_{1222}=t_{1333}=t_{2223}=-1$.   For $x=(\frac12,-\frac12,1)^\top$, we have
$$\aligned 
\mathcal{T}x^4=&(x_1+x_2+x_3)^4-8(x_1x_2^3+x_1x_3^3+x_2^3x_3)+6(x_1^2x_2^2+x_1^2x_3^2+x_2^2x_3^2)-24x_1x_2^2x_3\\\
=&1-8\left(-\dfrac1{2^4}+\dfrac12-\dfrac1{2^3}\right)+6\left(\dfrac1{2^4}+\dfrac1{2^2}+\dfrac1{2^2}\right)-3=-\dfrac98<0.
\endaligned$$

Case 3. $t_{1123}=t_{1223}=t_{1233}=1$ and  $t_{1112}=t_{2333}=t_{1113}=1$ and $ t_{1222}=t_{1333}=t_{2223}=-1$. Take $x=(\frac15,-\frac15,1)^\top$. Then we have $$ \aligned 
\mathcal{T}x^4=&(x_1+x_2+x_3)^4-8(x_1x_2^3+x_1x_3^3+x_2^3x_3)+6(x_1^2x_2^2+x_1^2x_3^2+x_2^2x_3^2)\\
=&1-8\left(-\dfrac1{5^4}+\dfrac15-\dfrac1{5^3}\right)+6\left(\dfrac1{5^4}+\dfrac1{5^2}+\dfrac1{5^2}\right)=-\dfrac{21}{625}<0.
\endaligned$$
The above three cases imply that the equality, $t_{1123}=t_{1223}=t_{1233}=-1,$ is necessary.

Case 4. $t_{1123}=t_{1223}=t_{1233}=-1,$ but $ t_{1222}=t_{2333}=t_{1113}$ and $ t_{1112}=t_{1333}=t_{2223}$ can't hold. Without loss the generality, let $t_{1112}=t_{2333}=t_{1113}=1$ and $ t_{1222}=t_{1333}=t_{2223}=-1$. Take $x=(-1,-3,-1)^\top$. Then we have $$ \aligned 
\mathcal{T}x^4=&(x_1+x_2-x_3)^4+8(x_1^3x_3+x_2x_3^3-x_1x_2^3)+6(x_1^2x_2^2+x_1^2x_3^2+x_2^2x_3^2)-24x_1x_2x_3^2\\
=&(-1-3+1)^4-8(1+3-3^3)+6(9+1+9)-24\times3=-61<0.
\endaligned$$
This is a contradiction to  the positive definiteness of $\mathcal{T}$, and hence, the  conditions (\textbf{III})  are  necessary.

Sufficiency.  $t_{1123}=t_{1223}=t_{1233}= -1$, $ t_{1222}=t_{2333}=t_{1113}$ and $ t_{1112}=t_{1333}=t_{2223}$.  Without loss the generality, let $t_{1222}=t_{2333}=t_{1113}=1$ and $ t_{1112}=t_{1333}=t_{2223}=-1$. 	Rewriting  $\mathcal{T}x^4 $ as follow, 
$$\aligned 
\mathcal{T}x^4=&(x_1+x_2-x_3)^4+8(x_1^3x_3+x_2x_3^3-x_1^3x_2)+6(x_1^2x_2^2+x_1^2x_3^2+x_2^2x_3^2)-24x_1x_2x_3^2\\
\ge&(x_1+x_2-x_3)^4+8(x_1^3x_3+x_2x_3^3-x_1^3x_2)+5(x_1^2x_2^2+x_1^2x_3^2+x_2^2x_3^2)-24x_1x_2x_3^2.
\endaligned$$
By (i), $\mathcal{T}x^4\ge0$ for all $x\in\mathbb{R}^3$.  It is not difficult to verify that the solutions of the equation $\mathcal{T}x^4=0$ is only original point $O(0,0,0).$  So,  $\mathcal{T}x^4>0$ for all $x\ne 0.$ That is, $\mathcal{T}$ is positive definite.

(iii)   Assume $t_{iijj} =2.5$ for all $i,j\in\{1,2,3\}$ and $i\ne j$.  If the conditions {\bf (IV)} can't hold when $\mathcal{T}$ is positive definite, then there are only two cases.

Case 1.  one of $\{t_{1123}, t_{1223}, t_{1233} \}$ is only $-1$. We might be $t_{1123}=t_{1233}=1$ and $ t_{1223}=-1$ and $t_{1112}=t_{2333}=t_{1113}=1$ and $ t_{1222}=t_{1333}=t_{2223}=-1$.  Then for $x=(\frac14,-\frac14,1)^\top$, we have
$$\aligned 
\mathcal{T}x^4=&(x_1+x_2+x_3)^4-8(x_1x_2^3+x_1x_3^3+x_2^3x_3)+9(x_1^2x_2^2+x_1^2x_3^2+x_2^2x_3^2)-24x_1x_2^2x_3\\
=&1-8\left(-\dfrac1{4^4}+\dfrac14-\dfrac1{4^3}\right)+9\left(\dfrac1{4^4}+\dfrac1{4^2}+\dfrac1{4^2}\right)-24\times \dfrac1{4^3}=-\dfrac{15}{256}<0.
\endaligned$$

Case 2.  $t_{1123}=t_{1223}=t_{1233}=-1,$ but $ t_{1222}=t_{2333}=t_{1113}$ and $ t_{1112}=t_{1333}=t_{2223}$ can't hold. Without loss the generality, let $t_{1112}=t_{2333}=t_{1113}=1$ and $ t_{1222}=t_{1333}=t_{2223}=-1$. Take $x=(-1,-3,-1)^\top$. Then we have $$ \aligned 
\mathcal{T}x^4=&(x_1+x_2-x_3)^4+8(x_1^3x_3+x_2x_3^3-x_1x_2^3)+9(x_1^2x_2^2+x_1^2x_3^2+x_2^2x_3^2)-24x_1x_2x_3^2\\
=&(-1-3+1)^4-8(1+3-3^3)+9(9+1+9)-24\times3=-4<0.
\endaligned$$
This obtains a contracdiction, and so, the conditions {\bf (IV)} is necessary.

Now we show the sufficiency. The second condition easily established by the proof of (ii), we only show the conclusion holds under the conditions: $t_{1123}=t_{1223}=t_{1233}= 1$ and $\mbox{two of } \{t_{1123}, t_{1223}, t_{1233}\} \mbox{ are  }-1$.  Let $t_{1222}=t_{2223}=t_{1113}=1$ and $ t_{1112}=t_{1333}=t_{2333}=-1$.

Condition: $t_{1123}=t_{1223}=t_{1233}= 1$. 	Rewriting  $\mathcal{T}x^4 $ as follow, 
$$ 
\mathcal{T}x^4=(x_1+x_2+x_3)^4-8(x_1^3x_2+x_1x_3^3+x_2x_3^3)+9(x_1^2x_2^2+x_1^2x_3^2+x_2^2x_3^2).
$$
It is easy to verify the global  minimum problem, $\min\mathcal{T}x^4$
has unique minimum $0$ at the origin coordinates $O(0,0,0)$. So,   $\mathcal{T}$ is positive definite. 	

Condition: $\mbox{two of } \{t_{1123}, t_{1223}, t_{1233}\} \mbox{ are  }-1$.  Without loss the generality, take $t_{1123}=t_{1223}=-1$ and $  t_{1233}=1. $ Then we have 
$$\mathcal{T}x^4=(x_1+x_2-x_3)^4+8(x_1^3x_3+x_2^3x_3-x_1^3x_2)+9(x_1^2x_2^2+x_1^2x_3^2+x_2^2x_3^2).$$
Solve the global  minimum problem, $\min\mathcal{T}x^4$, to yield its minimum value $0$ at the origin coordinates $O(0,0,0)$. So,   $\mathcal{T}$ is positive definite. 	

(iv)  Assume $t_{iijj} \ge\dfrac83$ for all $i,j\in\{1,2,3\}$ and $i\ne j$. Obvously,  the condition $|t_{iijk}|=1$ for all $i,j,k\in\{1,2,3\}, i\ne j, i\ne k, k\ne j,$ is equivalent to 
$$\begin{cases}
	t_{1123}=t_{1223}=t_{1233}=1; \mbox{ or}\\
	t_{1123}=t_{1223}=t_{1233}=-1;\mbox{ or}\\
	\mbox{two of } \{t_{1123}, t_{1223}, t_{1233}\} \mbox{ are  }-1;\mbox{ or}\\
	\mbox{one of } \{t_{1123}, t_{1223}, t_{1233}\} \mbox{ are  }-1.
\end{cases} \leqno{(\textbf{V})}$$ We need only show the conclusion holds under the conditions: one of $\{t_{1123}, t_{1223}, t_{1233} \}$ is $-1$ or $t_{1123}=t_{1223}=t_{1233}=-1$. Other two cases directly follow from (iii). Without loss the generality, let $t_{1222}=t_{2223}=t_{1113}=1$ and $ t_{1112}=t_{1333}=t_{2333}=-1$. 

Assume one of $\{t_{1123}, t_{1223}, t_{1233} \}$ is only $-1$. We might take $t_{1123}= t_{1223}=1$ and $ t_{1233}=-1$. 
Then we have 
$$\mathcal{T}x^4\geq(x_1+x_2+x_3)^4-8(x_1^3x_2+x_1x_3^3+x_2x_3^3)+10(x_1^2x_2^2+x_1^2x_3^2+x_2^2x_3^2)-24x_1x_2x_3^2.$$ Let
$$g(x_1,x_2,x_3)=(x_1+x_2+x_3)^4-8(x_1^3x_2+x_1x_3^3+x_2x_3^3)+10(x_1^2x_2^2+x_1^2x_3^2+x_2^2x_3^2)-24x_1x_2x_3^2.$$ Solve the global  minimum problem, $$\min \left\{g(x_1,x_2,x_3);x=(x_1,x_2,x_3)^\top\in\mathbb{R}^n\right\}$$ to yield its minimum value $0$ at the origin coordinates $O(0,0,0)$. So,  $\mathcal{T}x^4\geq g(x_1,x_2,x_3)>0$ for all $x\in\mathbb{R}^3\setminus\{0\}$. Thus $\mathcal{T}$ is positive definite. 	

Assume $t_{1123}=t_{1223}=t_{1233}=-1$.   Then we have 
$$\mathcal{T}x^4\geq(x_1+x_2-x_3)^4+8(x_1^3x_3+x_2x_3^3-x_1^3x_2)+10(x_1^2x_2^2+x_1^2x_3^2+x_2^2x_3^2)-24x_1x_2x_3^2.$$ Let
$$f(x_1,x_2,x_3)=(x_1+x_2-x_3)^4+8(x_1^3x_3+x_2x_3^3-x_1^3x_2)+10(x_1^2x_2^2+x_1^2x_3^2+x_2^2x_3^2)-24x_1x_2x_3^2.$$ Solve the global  minimum problem, $$\min \left\{f(x_1,x_2,x_3);x=(x_1,x_2,x_3)^\top\in\mathbb{R}^n\right\}$$ to yield its minimum value $0$ at the origin coordinates $O(0,0,0)$. So,  $\mathcal{T}x^4\geq f(x_1,x_2,x_3)>0$ for all $x\in\mathbb{R}^3\setminus\{0\}$. That is, $\mathcal{T}$ is positive definite.
This  completes the proof.	\\

By applying Theorems \ref{thm:31} (i) and (ii),  the following  inequalities are established easily  for  ternary quartic homogeneous polynomials.

\begin{corollary}\label{cor:32}	If $(x_1,x_2,x_3)\ne(0,0,0)$, then 
	\begin{itemize}
		\item [(i)] $(x_1+x_2-x_3)^4+5(x_1^2x_2^2+x_1^2x_3^2+x_2^2x_3^2)\geq8(x_1^3x_2-x_1^3x_3-x_2x_3^3)+24x_1x_2x_3^2,
		$\\ with equality if and only if  $x_1= x_2=x_3$;
		\item [(ii)] $(x_1+x_2-x_3)^4+6(x_1^2x_2^2+x_1^2x_3^2+x_2^2x_3^2)>8(x_1^3x_2-x_1^3x_3-x_2x_3^3)+24x_1x_2x_3^2.$
	\end{itemize} 
	Furthermore, these (strict) inequalities still hold if $x_1^3x_2$ and $x_1x_2^3$,  $x_1^3x_3$ and $x_1x_3^3$, $x_2x_3^3$ and  $x_2^3x_3$ are simultaneously  exchangeable.
\end{corollary}
By applying Theorems \ref{thm:31} (iii) and (iv),  the following strict inequalities are established easily  for  ternary quartic homogeneous polynomials.
\begin{corollary}\label{cor:33}	If $(x_1,x_2,x_3)\ne(0,0,0)$, then 
	\begin{itemize}
		\item [(i)] $(x_1+x_2+x_3)^4+9(x_1^2x_2^2+x_1^2x_3^2+x_2^2x_3^2)>8(x_1x_3^3+x_1^3x_2+x_2^3x_3);
		$ 
		\item [(ii)] $(x_1+x_2+x_3)^4+10(x_1^2x_2^2+x_1^2x_3^2+x_2^2x_3^2)>8(x_1x_3^3+x_1^3x_2+x_2^3x_3)+24x_1x_2x_3^2;
		$ \item [(iii)] $(x_1+x_2-x_3)^4+10(x_1^2x_2^2+x_1^2x_3^2+x_2^2x_3^2)>8(x_1^3x_2-x_1^3x_3-x_2x_3^3)+24x_1x_2x_3^2;$
		\item [(iv)] $(x_1+x_2-x_3)^4+9(x_1^2x_2^2+x_1^2x_3^2+x_2^2x_3^2)>8(x_1^3x_2-x_1^3x_3-x_2x_3^3).$
	\end{itemize} 
	Furthermore, these strict inequalities still hold if $x_1^3x_2$ and $x_1x_2^3$ are exchangeable, or  $x_1^3x_3$ and $x_1x_3^3$ are exchangeable, or $x_2x_3^3$ and  $x_2^3x_3$ are exchangeable.
\end{corollary}

\section{Conclusions}

For a 4th order 3-dimensional symmetric tensor with its entries  $1$ or $-1$,   the analytic  necessary and sufficient conditions are established for its positive definiteness.   Several  (strict) inequalities  of ternary quartic homogeneous polynomial are built by means of these analytic conditions.


\begin{thebibliography}{00}
	\bibitem{K2016}K. Kannike,  Vacuum stability of a general scalar potential of a few fields. Eur. Phys. J. C. 76(2016), 324; Erratum. Eur. Phys. J. C. 78(2018), 355
	\bibitem{SQ2024}Y. Song, L. Qi,  Boundedness from below conditions for a general scalar potential of two real scalars fields and the Higgs boson, to appear:  Theor. Math. Phys.  chinaXiv:202011.00118, arXiv: 2011.11262
	\bibitem{Q2005}  L. Qi, Eigenvalues of a real supersymmetric tensor, J. Symbolic Comput. 40(2005) 1302-1324
	\bibitem{R1922}E. L. Rees,  Graphical Discussion of the Roots of a Quartic Equation, Amer. Math. Monthly,  {\bf 29(2)}, (1922) 51-55
	\bibitem{L1988}D.  Lazard, Quantifier elimination: Optimal solution for two classical examples, J.  Sym. Comput.	{\bf 5(1-2)} (1988) 261-266 
	\bibitem{GL1964} R.N. Gadem, C.C.  Li,   On positive definiteness of quartic forms of twovariables,  IEEE Trans. Autom. Control.   {\bf AC-9}, (1964) 187-188
	\bibitem{K1965}  W.H. Ku,  Explicit criterion for the positive definiteness of a general quartic form,  IEEE Tram. Autom. Control. {\bf AC-10(3)}, (1965) 372-373   
	\bibitem{JM1981} E.I. Jury, M.   Mansour, Positivity and nonnegativity of a quartic equation and related problems, IEEE Trans. Autom. Control.  {\bf AC-26(2)}, (1981) 444-451 
	\bibitem{WQ2005}F. Wang,  L. Qi,  Comments on ``Explicit criterion for the positive definiteness of a general quartic form",  IEEE Trans. Autom. Control.  {\bf 50(3)}, (2005) 416-418
	\bibitem{HH1996} M. A. Hasan,   A. A.  Hasan,   A procedure for the positive definiteness of forms of even order.  IEEE Trans. Auto. Control.   {\bf 1996} {\em 41(4)}, 615-617.
	\bibitem{F1998}M. Fu,  Comments on ``A procedure for the positive definiteness of forms of even order". IEEE Trans. Auto. Control.  {\bf 1998} {\em 43(10)},  1430.
	
	\bibitem{G2021}Y.  Guo, A necessary and sufficient condition for the positive definite problem of a binary quartic form. J. Math.  {\bf 2021} {\em 2021}, 2339746.
	
	\bibitem{QSZ2022}L. Qi, Y. Song,  X. Zhang, Positivity Conditions for Cubic, Quartic and Quintic Polynomials,  J. Nonlinear Convex Anal.   {\bf 23(2)}, (2022) 191-213
	\bibitem{SQ2021} Y. Song, L. Qi,  Analytical expressions of copositivity for 4th order symmetric tensors, Analy. Appl.   {\bf 19(5)}, (2021)  779-800.
	\bibitem{NQW2008}Q. Ni,   L. Qi, F.  Wang,  An eigenvalue method for testing the positive definiteness of a multivariate form. IEEE Trans Automatic Control, 2008, 53: 1096-1107
	\bibitem{NQZ2009}  M.  Ng, L. Qi, G.  Zhou, Finding the largest eigenvalue of a non-negative tensor, SIAM J Matrix Anal Appl, 2009, 31: 1090-1099
	\bibitem{S2021}Y. Song,  Positive definiteness for 4th order symmetric tensors and applications, Analy. Math.  Physics, 10 (2021), 11
		\bibitem{UW1994} G. Ulrich,  L.T.  Watson,   Positivity conditions for quartic polynomials. SIAM J. Sci. Comput. 15(1994), 528-544
\end{thebibliography}
\end{document}